\documentclass[14pt, a4paper]{amsart}

\usepackage[14pt]{extsizes}
\usepackage{amsmath}
\usepackage{amssymb}
\usepackage{enumerate}

\newtheorem{Theorem}{Theorem}

\begin{document}

\author{A.\,N.~Lavrov}
\address{Higher School of Modern Mathematics, MIPT, Moscow, Klimentovsky per. 1}
\email{lavrov.an@mipt.ru}

\title{Einstein homogeneous spaces associated \\ with reducible Clifford modules}

\maketitle

\begin{abstract}
    V. C\'{o}rtes suggested the construction of noncompact pseudo-Riemannian homogeneous spaces which naturally carry
    an invariant quaternion-pseudo-K\"{a}hler structure. The construction is based on consideration of pseudo-Euclidean
    space $V = \mathbb{R}^{p,q}$, $\text{Cl}^{0}(V)$-module $W$, $V$-valued $\mathfrak{so}(V)$-equivariant skew-symmetric
    bilinear form $\Pi$ on $W$ and pseudo-Riemannian metric associated with the orthogonal decomposition
    $V = \mathbb{R}^{3, 0} \oplus \mathbb{R}^{p-3, q}$. The corresponding metric is Riemannian if and only if $p = 3$.
    Moreover, any noncompact quaternion-K\"{a}hler homogeneous space can be obtain by this construction for $p = 3$.
    Earlier D. V. Alekseevsky and the author generalized the construction of V. C\'{o}rtes for $p = 3$ to the case
    $p \equiv 3 \pmod{4}$, irreducible $\text{Cl}^{0}(V)$-module $W$ and Riemannian metric associated with the decomposition
    $V = \mathbb{R}^{p, 0} \oplus \mathbb{R}^{0, q}$. It was shown that doing in this way we obtain the series of noncompact
    Riemannian Einstein homogeneous spaces $M_{p, q}$. The current paper is devoted to further generalization of this result
    to the case $p \equiv 3 \pmod{4}$, reducible $\text{Cl}^{0}(V)$-module $W$ and a family of pseudo-Riemannian
    metrics associated with the decomposition $V = \mathbb{R}^{p, 0} \oplus \mathbb{R}^{0, q}$. As a result we extend
    the series $M_{p, q}$ to larger series of noncompact pseudo-Riemannian Einstein homogeneous spaces $M_{p, q}(m)$
    where $m$ is the tuple of four (eight) non-negative integers for $q \not\equiv 3 \pmod{4}$
    (for $q \equiv 3 \pmod{4}$, respectively). 

\vspace{2mm}
\end{abstract}

\footnotetext{This work is supported by the Russian Science Foundation \textnumero 25-21-00918,
    https://rscf.ru/project/25-21-00918/}
    
\section{Introduction}

    The current paper is devoted to the construction of the infinite series of pseudo-Riemannian Einstein homogeneous spaces
    with negative scalar curvature. The construction is based on the construction considered in \cite{C} of noncompact
    quaternion-pseudo-K\"{a}hler homogeneous spaces which is the following. For pseudo-Euclidean space $V = \mathbb{R}^{p, q}$,
    $p \geq 3$, consider a specific $\mathbb{Z}_{2}$-graded extension $\mathfrak{p}_{p, q}(W, \Pi)$ of Poincar\'{e} algebra
    $\text{Lie}\big( \text{Isom}(V) \big)$, where $\text{Isom}(V)$ is the isometry group of $V$, by arbitrary module $W$
    over the even Clifford algebra $\text{Cl}^{0}(V)$ and $V$-valued $\mathfrak{so}(V)$-equivariant skew-symmetric
    bilinear form $\Pi: \wedge^{2} W \longrightarrow V$. The resulting Lie algebra $\mathfrak{p}_{p, q}(W, \Pi)$
    is called an \textit{extended Poincar\'{e} algebra}, all Lie algebras of this form for any signature $(p, q)$ were classified
    in \cite{AC}. Next, consider the extension $\mathfrak{g}_{p, q}(W, \Pi)$ of Lie algebras $\mathfrak{p}_{p, q}(W, \Pi)$
    by the derivation $D$ with eigenvalues $(0, \frac{1}{2}, 1)$. Let $G_{p, q}(W, \Pi)$ be the corresponding simply
    connected Lie group, while $H \subset \text{Spin}^{0}(V) \subset G(W, \Pi)$ be its Lie subgroup which respects
    the orthogonal decomposition $V = \mathbb{R}^{3, 0} \oplus \mathbb{R}^{p-3, q}$. Then the quotient manifold $G(W, \Pi)/H$
    naturally admits an invariant quaternion-pseudo-K\"{a}hler structure. This structure is quaternion-K\"{a}hler if and
    only if $p = 3$. Moreover, any noncompact quaternion-K\"{a}hler homogeneous space can be obtained by this construction
    for $p = 3$.

    In \cite{AL} the construction of V. C\'{o}rtes for $p = 3$ is generalized to the case $p \equiv 3 \pmod{4}$. The main point of
    generalization is consideration of specific Riemannian metric associated with the orthogonal decomposition $V = \mathbb{R}^{p, 0} \oplus \mathbb{R}^{0, q}$
    instead of $V = \mathbb{R}^{3, 0} \oplus \mathbb{R}^{p-3, q}$ as in the original construction. However, in \cite{AL}
    the case of irreducible $\text{Cl}^{0}(V)$-module is only considered. As a result the new series of noncompact Riemannian
    Einstein homogeneous spaces $M_{p,q}$ which are no longer quaternion-K\"{a}hler is constructed in \cite{AL}.
    The current paper is devoted to further generalization of this result, namely, we consider the same construction as in
     \cite{AL}, but $\text{Cl}^{0}(V)$-module is supposed to be reducible and instead of a single metric we consider some family
     of pseudo-Riemannian metrics on the corresponding homogeneous space. For the case $p \equiv 3 \pmod{4}$,
     $q \not\equiv 3 \pmod{4}$ the considered family of pseudo-Riemannian homogeneous spaces can be parametrized
     by a pair $(A, \hat{A})$ of nondegenerate symmetric matrices of the same size. The main results are the following:
    \begin{enumerate}
        \item The pseudo-Riemannian homogeneous space associated with the pair $(A, \hat{A})$ is an Einstein manifold
            if and only if the matrix equality $( \hat{A}^{-1} A \big)^{2} = \text{Id}$ is satisfied (see Theorem \ref{EinsteinCriterion}).
        \item The Einstein homogeneous space associated with the pair $(A, \hat{A})$ is isometric to the homogeneous space associated
            with the pair of diagonal matrices with elements $\pm 1$ on the diagonal (see Theorem \ref{DiagonalizationTheorem}).
    \end{enumerate}
    In particular, the last result implies that the constructed Einstein homogeneous spaces can be parametrized up to isometry
    by the tuple of four non-negative integers (see Theorem \ref{TupleParametrization}). For the case $p, q \equiv 3 \pmod{4}$
    the considered family is parametrized by two pairs $(A^{+}, \hat{A}^{+})$, $(A^{-}, \hat{A}^{-})$ of nondegenerate
    symmetric matrices. The similar results are also valid (see Theorems \ref{EinsteinCriterion_q3}, \ref{DiagonalizationTheorem_q3})
    and the corresponding Einstein homogeneous spaces can be parametrized by the tuple of eight non-negative integers
    (see Theorem \ref{TupleParametrization_q3}).

\section{Construction of homogeneous space}
  
    Let $V = \mathbb{R}^{p, q}$ be the standard pseudo-Euclidean vector space with indefinite scalar product $\langle \cdot, \cdot \rangle$
    of the signature $(p, q)$. We will denote the dimension of $V$ by $n := \dim V = p + q$. Remind that the Clifford algebra
    $\text{Cl}(V)$ is the quotient of the tensor algebra $T(V)$ by the two-sided ideal\footnote{Note that the other
    definition of Clifford algebra is often used in the literature in which the considered ideal is generated by elements 
    $v \otimes v - \langle v, v \rangle \cdot 1$, $v \in V$.} generated by the elements of the form $v \otimes v + \langle v, v \rangle \cdot 1$, $v \in V$.
    The involution $v \mapsto (-v)$ on $V$ induces the $\mathbb{Z}_{2}$-grading $\text{Cl}(V) = \text{Cl}^{0}(V) + \text{Cl}^{1}(V)$
    where $\text{Cl}^{0}(V)$ is an even subalgebra often called the even Clifford algebra. The group $\text{Pin}(V) \subset \text{Cl}(V)$
    generated by all vectors $v \in V$ with $\langle v, v \rangle = \pm 1$ is called $\text{Pin}$-group. Its subgroup
    $\text{Spin}(V) := \text{Pin}(V) \cap \text{Cl}^{0}(V)$ consisted of even elements is called spinor group.
    
    The Lie algebra $\mathfrak{spin}(V) := \text{Lie}\big( \text{Spin}(V) \big) \simeq \text{span}\{ [u, v] \ | \ u, v \in V \} \subset \text{Cl}^{0}(V)$, 
    where $[u, v] := uv - vu$ is the commutator inside the Clifford algebra is isomorphic to the Lie algebra $\mathfrak{so}(V) := \text{Lie}\big( SO(V) \big)$
    which can be identified with the space of bivectors $\wedge^{2} V$ by the formula 
    \begin{equation}
        u\wedge v : w \mapsto  \langle v, w \rangle u -   \langle u, w  \rangle v. 
    \end{equation}
    The isomorphism $\mathfrak{spin}(V) \simeq \wedge^{2} V$ is given by the map $[u, v] \mapsto (-4) \cdot (u \wedge v)$
    for any $u, v \in V$. In particular, $u \cdot v \mapsto (-2) \cdot (u \wedge v)$ if $u \perp v \in V$.

    Next, consider the orthogonal decomposition $V = \mathbb{R}^{p, q} = E \oplus E'$ where $E = \mathbb{R}^{p, 0}$ and $E' = \mathbb{R}^{0, q}$.
    Let $(e_{1}, .., e_{p},$ $e'_{1}, .., e'_{q})$ be the standard basis of $V$ such that $(e_{1}, .., e_{p})$ is the standard basis of $W$,
    while $(e'_{1}, .., e'_{q})$ is the standard basis of $E' := \mathbb{R}^{0, q} \subset V$. Note that the considered decomposition
    induces the decomposition of the Lie algebra $\mathfrak{so}(V) \simeq V \wedge V = E \wedge E + E' \wedge E' + E \wedge E'$.
    Also note that the basis $(e_{1}, .., e_{p},$ $e'_{1}, .., e'_{q})$ determines the orientation of the vector spaces $V$, $E$ and $E'$.
    Let $\hat{K} \subset \text{Spin}(V)$ be the subgroup of $\text{Spin}(V)$ consisted of all elements which respects
    the decomposition $V = E \oplus E'$. One can show that $\hat{K}$ is a maximal compact subgroup of $\text{Spin}(V)$.

    Starting from this point to Section 4 we will consider only the case $p \equiv 3 \pmod{4}$, $q \not\equiv 3 \pmod{4}$.
    One can show that for this case there exists only one irreducible $\text{Cl}^{0}(V)$-module $\Sigma$ of the dimension
    $d := \dim{\Sigma}$. So any $\text{Cl}^{0}(V)$-module $W$ can be represented as $W \simeq \mathbb{R}^{l} \otimes \Sigma$.
    Denote the standard basis of $\mathbb{R}^{l}$ by $(\xi_1, ..., \xi_l)$. According to \cite[Lemma 3]{C},
    $\Sigma$ is an irreducible representation of the compact subgroup $\hat{K} \subset \text{Spin}(V) \subset \text{Cl}^{0}(V)$.
    Consequently, there exists an unique (up to multiplication by a positive scalar) $\text{Ad}_{\hat{K}}$-invariant symmetric
    positive-definite bilinear form $b$ on $\Sigma$. So any $\text{Ad}_{\hat{K}}$-invariant symmetric bilinear form on $W$
    has the form $b_{A} = \epsilon_{A} \otimes b$ where $\epsilon_{A}$ is some symmetric bilinear form on $\mathbb{R}^{l}$
    defined by the matrix $A$ in the basis $(\xi_{k})$. Further we will assume that $\epsilon_{A}$ is non-degenerated..

    Using the bilinear form $b_{A}$ we define the $\mathfrak{so}(V)$-equivariant linear map $\Pi_{A}:\bigwedge^2 W \longrightarrow V$
    by the following formula
    \begin{equation}\label{V-valuedForm}
        \langle \Pi_{A}(s, s'), v \rangle = b_{A}( \omega v s, s')
    \end{equation}
    which defines $\mathbb{Z}_{2}$-graded Lie algebra 
    \begin{equation}\label{ExtendedPoincareAlgebra}
        \mathfrak{p}_{p, q}(A) := \mathfrak{p}_{0} + \mathfrak{p}_{1} = \big( \mathfrak{so}(V) + V \big) + W.
    \end{equation} 

    Lie algebra $\mathfrak{p}_{p, q}(A)$ admits the derivation $D$ such that
    \begin{equation}\label{Derivation}
        \text{ad}_{D}|_{V} = \text{Id}, \ \ \ \text{ad}_{D}|_{W} = \frac{1}{2}\text{Id}, \ \ \ \text{ad}_{D}|_{\mathfrak{so}(V)} = 0.
    \end{equation} 
    Using the derivation $D$ we can define $\mathbb{Z}_{2}$-graded Li algebra 
    \begin{equation}\label{ExtensionOfEPA}
        \mathfrak{g}_{p, q}(A) := \mathfrak{g}_{0} + \mathfrak{g}_{1} 
            = \big( \mathbb{R}D + \mathfrak{p}_{0} \big) + \mathfrak{p}_{1}
            = \big( \mathbb{R}D + \mathfrak{so}(V) + V \big) + W.
    \end{equation}
    Let $G_{p,q}(A)$ be the corresponding simply connected Lie group. Note that the identity component of the spinor group
    $\text{Spin}_{0}(V) \subset \text{Spin}(V)$ can be considered as a closed Lie subgroup of $G_{p, q}(A)$ corresponding
    to Lie subalgebra $\mathfrak{so}(V) \subset \mathfrak{g}_{p, q}(A)$ (see \cite[Cor. 3]{C}).

    Denote the maximal connected subgroup of $\text{Spin}_{0}(V)$ which respects the orthogonal decomposition $V = E \oplus E'$
    by $K \subset \text{Spin}_{0}(V)$. We have that $K = \text{Spin}(E) \times \text{Spin}(E') \simeq \text{Spin}(p) \times \text{Spin}(q)$,
    in particular, $K = \hat{K} \cap \text{Spin}_{0}(V)$ is a maximal compact subgroup of $\text{Spin}_{0}(V)$ and it corresponds
    to Lie subalgebra $\mathfrak{k} := \mathfrak{so}(E) + \mathfrak{so}(E') \simeq E \wedge E + E' \wedge E' \subset \mathfrak{so}(V)$.
    As well as the subgroup $\text{Spin}_{0}(V) \subset G_{p, q}(A)$, the group $K$ can be considered as a closed Lie subgroup of $G_{p, q}(A)$.
    Now consider the homogeneous space $M_{p,q}(A) := G_{p, q}(A)/K$ with the reductive decomposition
    \begin{equation}\label{ReductiveDecomposition}
        \mathfrak{g}_{p, q}(A) = \mathfrak{k} + \mathfrak{m}, \ \ \text{where} \ \mathfrak{k} = E \wedge E + E' \wedge E',
            \ \mathfrak{m} = \mathbb{R}D + E \wedge E' + V + W.
    \end{equation}
    Note that $\mathbb{RD}$, $E \wedge E'$, $E$, $E'$ are irreducible $\text{Ad}_{K}$-modules. Moreover,
    $\mathbb{RD}$, $E \wedge E'$, $E$, $E'$ and $W$ are pairwise non-equivalent for $q > 1$. For the case $q = 1$
    we have the isomorphisms of $\text{Ad}_{K}$-modules $E' \simeq \mathbb{R}D$ and $E \wedge E' \simeq E$.
    Next, the tangent space $M_{p, q}(A)$ at the point $o := [eK] \in M_{p, q}(A)$ is isomorphic to the subspace
    $\mathfrak{m} \subset \mathfrak{g}_{p, q}(A)$, i. e. $T_{o} M_{p, q}(A) \simeq \mathfrak{m}$.
    The $G_{p,q}(A)$-invariant pseudo-Riemannian metrics on $M_{p, q}(A)$ are in one-to-one correspondence with
    the $\text{Ad}_{K}$-invariant scalar products on $\mathfrak{m}$. The Ricci curvature of such invariant metric
    at the point $o \in M_{p, q}(A)$ can be expressed in terms of the corresponding scalar product on $\mathfrak{m}$
    and the Lie bracket. 

    Now consider the positive-definite scalar product $g_{V}$ on $V$ such that $g_{V}|_{E} = \langle -, - \rangle|_{E}$,
    but $g_{V}|_{E'} = (-1)\cdot\langle -, - \rangle|_{E'}$. It induces the positive-definite scalar product $g_{E \wedge E'}$ on $E \wedge E'$. 
    Also consider the positive-definite metric $g_{\mathbb{RD}}$ on $\mathbb{R}D$ determined by the relation $g_{\mathbb{R}D}(D, D) = 1$.
    Now define the indefinite scalar product on the whole subspace $\mathfrak{m}$ by the following formula
    \begin{equation}\label{ScalarProduct}
        g_{\hat{A}} = g_{\mathbb{R}D} \oplus g_{E \wedge E'} \oplus g_{V} \oplus b_{\hat{A}},
    \end{equation}
    which depends on the arbitrary symmetric non-degenerate $l \times l$-matrix $\hat{A}$.
    Since $g_{\hat{A}}$ is $\text{Ad}_{K}$-invariant it can be extended to the pseudo-Riemannian metric $g_{M, \hat{A}}$
    on the homogeneous space $M_{p, q}(A) = G_{p, q}(A)/K$.

\section{Einstein equation}

    Ricci curvature can be expressed through so-called Nomizu operators. By definition, \textit{Nomizu operator} $L_{X}$
    associated with vector field $X$ on $M_{p, q}(A)$ is the tensor field given by the formula $L_{X}Y := - \nabla_{Y}X$ where
    $Y$ is arbitrary vector field. In particular, $L_{X}|_{T_{o}M_{p,q}(A)}$ is an endomorphism of the tangent space $T_{o}M_{p,q}(A) \simeq \mathfrak{m}$. 
    On the other hand, any element of the Lie algbera $x \in\mathfrak{g}_{p, q}(A)$ induces the Killing vector field $\alpha(x)$
    on $M_{p, q}(A)$ which satisfies $\alpha(x)|_{o} = \pi(x) \in \mathfrak{m} \subset \mathfrak{g}(W, \Pi)$ where
    $\pi: \mathfrak{g}_{p,q}(A) \longrightarrow \mathfrak{m}$ is the projection onto $\mathfrak{m}$ with respect to
    the reductive decomposition (\ref{ReductiveDecomposition}). For any element $x \in \mathfrak{g}_{p,q}(A)$ we denote
    the endomorphism  $L_{\alpha(x)}|_{T_{o}M_{p,q}}(A) \in \text{End}(\mathfrak{m})$ just by $L_{x}$. For simplicity we will also call
    such endomorphisms $L_{x}, x \in \mathfrak{g}_{p,q}(A),$ Nomizu operators. The endomorphism $L_{x}, x \in \mathfrak{g}_{p,q}(A),$
    satisfies the following equality for any $y, z \in \mathfrak{m}$ (see \cite{AVL})
    \begin{equation}\label{NomizuOperator}
        g(L_{x} y, z) = -\frac{1}{2} \Big( g\big(\pi([x, y]), z\big) - g\big(\pi([y, z]), x\big) + g\big(\pi([z, x]), y\big) \Big).
    \end{equation}
    From this one can see that Nomizu operators are skew-symmetric with respect to $g$. Also they satisfy the following equality
    \begin{equation}\label{NomizuIdentity}
        L_{y}\pi(x) - L_{x}\pi(y) = \pi([x, y]) \text{ for any } x, y \in \mathfrak{g}(\Pi).
    \end{equation}
    Along with Nomizu operators we will also use the operators $R_{x}: \mathfrak{g}_{p,q}(A)\rightarrow\mathfrak{m}, x \in \mathfrak{m},$
    defined as follows
    \begin{equation}
        R_{x} y := L_{y} x = L_{x}\pi(y) + \pi([x, y]). 
    \end{equation}
    The Riemann curvature tensor of the metric $g_{M_{p,q}(A)}$ at $o \in M_{p,q}(A)$ can be expressed through  Nomizu operators
    like this (see \cite[Sec. 2.3]{C})
    \begin{equation}\label{RiemannTensorFormula}
        R(x, y) = [L_{x}, L_{y}] + L_{[x, y]} \text{ for any } x, y \in \mathfrak{m}.
    \end{equation}

    Also consider the endomorphism $C \in \text{End}(\mathbb{R}^{l})$ which transforms the bilinear form $\epsilon_{\hat{A}}$
    to$\epsilon_{A}$, i. e. satisfies the equality
    \begin{equation}\label{BilinearFormChange}
        \epsilon_{\hat{A}}(C \xi, \xi' ) = \epsilon_{A}(\xi, \xi') \text{ for any } \xi, \xi' \in \mathbb{R}^{l}.
    \end{equation}
    From this equality it follows that $C$ is symmetric with respect to both bilinear forms $\epsilon_{\hat{A}}$ and $\epsilon_{A}$,
    i. e. the following equalities hold
    \begin{equation}\label{SymmetryOfEndomorphism}
        \epsilon_{\hat{A}}(C \xi, \xi' ) = \epsilon_{\hat{A}}(\xi, C \xi' ),  \ \ \ \epsilon_{A}(C \xi, \xi' ) = \epsilon_{A}(\xi, C \xi' ).
    \end{equation}
    In the basis $(\xi_{k})$ the endomorphism $C$ is given by the matrix $\hat{A}^{-1} A$. From (\ref{ScalarProduct}) we have that
    \begin{equation*}
        g_{\hat{A}}\Big( (C \otimes \text{Id}) s, s' \Big) = b_{\hat{A}}(Cs, s') = b_{A}(s, s').
    \end{equation*}

    \begin{Theorem}\label{NomizuOperators}
        For the metric $g_{M, \hat{A}}$ the Nomizu operators have the following form:
        \begin{enumerate}[(i)]
            \item $L_{x} = 0$ for any $x \in \mathbb{R}D + E \wedge E'$,
            \item $L_{x} = - \text{\emph{ad}}_{x}$ for any $x \in E \wedge E + E' \wedge E'$,
            \item $L_{e_{i}} = -D \otimes e_{i}^{*} + \sum\limits_{j=1}^{q} (e_{i} \wedge e'_{j}) \otimes {e'_{j}}^{*} 
                - \text{\emph{ad}}_{e_{i}} + \frac{1}{2} C \otimes \omega e_{i} $,
                where $1 \leq i \leq p$,
            \item $L_{e'_{j}} = -D \otimes {e'_{j}}^{*} + \sum\limits_{k=1}^{p} (e_{k} \wedge e'_{j}) \otimes e_{k}^{*} 
                - \text{\emph{ad}}_{e'_{j}} - \frac{1}{2} C \otimes \omega e'_{j} $,
                where $1 \leq j \leq q$,
            \item $L_{s} = -\frac{1}{2} D \otimes s^{*} 
                + \frac{1}{2} \sum\limits_{i=1}^{p}\sum\limits_{j=1}^{q} (e_{i} \wedge e'_{j}) \otimes \big( (e_{i} e'_{j}) s \big)^{*}
                - \frac{1}{2}\text{\emph{ad}}_{s}|_{W} - \text{\emph{ad}}_{s}|_{\mathbb{R}D + E \wedge E'}
                + \frac{1}{2}\sum\limits_{i=1}^{p} \Bigg( \Big( C \otimes \omega e_{i} \Big)  s \Bigg)  \otimes e_{i}^{*}
                - \frac{1}{2}\sum\limits_{j=1}^{q} \Bigg( \Big( C \otimes \omega e'_{j} \Big) s  \Bigg)   \otimes {e'_{j}}^{*}$ for any $s \in W$,
        \end{enumerate}
        where the identification $\mathfrak{m} \simeq \mathfrak{m}^{*}$ is given by the scalar product $g_{\hat{A}}$.
    \end{Theorem}
    \textit{Proof:}
        To compute $g(L_{x} y, z)$ we will use the formula (\ref{NomizuOperator}). Due to linearity it is enough to consider only
        the basis elements of $\mathfrak{m}$. Moreover, the form of the formula for $g(L_{x} y, z)$ essentially depends only
        on the pair of $\text{Ad}_{K}$-submodules $(F_{1}, F_{2})$, $F_{1}, F_{2} \in \{\mathbb{R}D, E, E', E \wedge E', W \}$,
        that the elements $y, z$ belong to.

        Note that (i), (ii) can be derived in the same way as (i), (ii) from \cite[Th. 4]{AL}. All necessary computations for
        $(L_{e_{i}} y, z)$ and $(L_{e'_{j}} y, z)$ for all pairs of $\text{Ad}_{K}$-submodules $(F_1, F_2)$ except $(W, W)$ also
        essentially coincide with the corresponding computations from \cite[Th. 4]{AL}. The computation for $(W, W)$
        is given below
        \begin{equation*}
            g_{\hat{A}}(L_{e_{i}} s, s') = \frac{1}{2} \langle [s, s'], e_{i} \rangle = \frac{1}{2} b_{A}(\omega e_{i} s, s')        
                = g_{\hat{A}} \Big( \frac{1}{2} (C \otimes \omega e_{i} ) s, s' \Big),
        \end{equation*}
        \begin{equation*}
            g_{\hat{A}}(L_{e'_{j}} s, s') = \frac{1}{2} \langle [s, s'], e'_{j} \rangle = \frac{1}{2} b_{A}(\omega e'_{j} s, s')        
                = g_{\hat{A}} \Big( \frac{1}{2} (C \otimes \omega e'_{j} ) s, s' \Big).
        \end{equation*}
        From these formulas we obtain the last summands in the expressions (iii) and (iv). For $(L_{s} y, z)$ the corresponding
        computations differ from those in \cite[Th. 4]{AL} by only pairs $(E, W)$ and $(E', W)$. For these pairs we have
        the formulas
        \begin{equation*}
            g_{\hat{A}}(L_{s} e_{i}, s') = \frac{1}{2}g_{\hat{A}}([s, s'], e_{i}) = \frac{1}{2} b_{A}(\omega e_{i} s, s')
                = g_{\hat{A}}\Bigg( \frac{1}{2} \Big( C \otimes \omega e_{i} \Big) s, s' \Bigg),
        \end{equation*}
        \begin{equation*}
            g_{\hat{A}}(L_{s}e'_{j}, s') = \frac{1}{2}g_{\hat{A}}([s, s'], e'_{j}) = - \frac{1}{2}b_{A}(\omega e'_{j} s, s')
                = g_{\hat{A}} \Big( - \frac{1}{2} \Big( C \otimes \omega e'_{j} \Big) s, s' \Big).
        \end{equation*}
        From these we obtain two last summands in the expression (v). 
    \hfill$\Box$

    \begin{Theorem}\label{EinsteinCriterion}
        Pseudo-Riemannian homogeneous space $\Big( M_{p, q}(A), g_{M, \hat{A}} \Big)$ is Einstein if and only if $C^{2} = \text{Id}$.
        The Einstein constant is $\lambda = - n  - \frac{1}{4} l d$.
    \end{Theorem}
    \textit{Proof:}
    From the formula (\ref{RiemannTensorFormula}) one can derive the following expression
    for the Ricci curvature at the point $o \in M_{p,q}(A)$ using Nomizu operators
    \begin{equation}\label{RicciQuadraticForm}
        \text{Ric}(x, x) = \text{tr}\big( - R_{x}^{2} + R_{L_{x}x} - R_{x}\text{ad}_{x}^{\mathfrak{k}} \big)
    \end{equation}
    where $\text{ad}_{x}^{\mathfrak{k}} := \pi^{\mathfrak{k}} \circ \text{ad}_{x}$, $\pi^{\mathfrak{k}}$ is the projection operator
    from $\mathfrak{g}_{p,q}(A)$ onto $\mathfrak{k}$ along $\mathfrak{m}$, i. e. $\pi^{\mathfrak{k}} = \text{Id} - \pi$.
    In the formula (\ref{RicciQuadraticForm}) and in all formulas below we consider the trace which is taken over the subspace
    $\mathfrak{m} \subset \mathfrak{g}_{p,q}(A)$.
        
    As in the previous theorem many computations which are necessary to prove the statement are completely similar
    to those from \cite[Th. 5]{AL}. In particular, from $L_{D} = 0$, $\text{ad}^{\mathfrak{k}}_{D} = 0$, and the equality
    $R_{x} = L_{x} \circ \pi + \pi \circ \text{ad}_{x}$, we obtain that
    \begin{equation*}
        \text{Ric}(D,D) = -\text{tr }R_{D}^{2} = -\text{tr} (\text{ad}_{D})^{2} = - \bigg( n + \frac{ld}{4} \bigg).
    \end{equation*}
    So if $\Big( M_{p, q}(A), g_{M, \hat{A}} \Big)$ is Einstein then the Einstein constant should be $\lambda = - \bigg( n + \frac{ld}{4} \bigg)$.

    Next, similarly to \cite[Th. 5]{AL} one can show that
    \begin{equation*}
        \text{tr } \pi \circ \text{ad}_{e_{i} \wedge e'_{j}} = 0,  \ \ \text{tr } \big( \pi \circ \text{ad}_{e_{i} \wedge e'_{j}} \big)^2 = 2 + \frac{l \cdot d}{4},
            \ \ \text{tr } R_{e_{i} \wedge e'_{j}}\text{ad}_{e_{i} \wedge e'_{j}}^{\mathfrak{k}} = n-2.
    \end{equation*}
    This implies that $\text{Ric}(e_{i} \wedge e'_{j}, e_{i} \wedge e'_{j}) = - \big( n + \frac{l \cdot d}{4} \big)$.
    For the elements of the form $L_{x} x$ the following formulas hold
    \begin{equation*}
        L_{e_{i}}e_{i} = L_{e'_{j}}e'_{j} = -D, \ \
            L_{s}s  = -\frac{1}{2} g_{\hat{A}}(s, s) D + \frac{1}{2} \sum\limits_{i,j=1}^{p,q} g_{\hat{A}}(e_{i}e'_{j}s, s) e_{i} \wedge e'_{j}.
    \end{equation*}
    Taking into account that $\text{ad}_{e_{i}}^{\mathfrak{k}} = \text{ad}_{e'_{j}}^{\mathfrak{k}} = \text{ad}_{s}^{\mathfrak{k}} = 0$
    we obtain
   \begin{equation}\label{Ric_E_E'}
        \text{Ric}(e_{i}, e_{i}) = -\text{tr } R_{e_{i}}^{2} - \Big( n + \frac{ld}{2} \Big) , \ \ \
        \text{Ric}(e'_{j}, e'_{j}) = -\text{tr } R_{e'_{j}}^{2} - \Big( n + \frac{ld}{2} \Big),
    \end{equation}
    \begin{equation}\label{Ric_S}
        \text{Ric}(s, s) = - \text{tr } R_{s}^{2} - \bigg( \frac{n}{2} + \frac{ld}{4} \bigg) g_{\hat{A}}(s, s).
    \end{equation}
 
    The computations of summands $\text{tr } R_{e_{i}}^{2}$, $\text{tr } R_{e_{i}}^{2}$, $\text{tr } R_{s}^{2}$ differ from those
    in \cite[Th. 5]{AL} so we present them here. First of all, note that
    \begin{equation*}
        R_{e_{i}} = L_{e_{i}} \circ \pi + \pi\circ\text{ad}_{e_{i}} = -D \otimes e_{i}^{*}
                + \sum\limits_{j=1}^{q} (e_{i} \wedge e'_{j}) \otimes {e'_{j}}^{*} + \frac{1}{2} C \otimes \omega e_{i},
    \end{equation*}
    \begin{equation*}
        R_{e'_{j}} = L_{e'_{j}} \circ \pi + \pi\circ\text{ad}_{e'_{j}} = -D \otimes {e'_{j}}^{*}
                + \sum\limits_{k=1}^{p} (e_{k} \wedge e'_{j}) \otimes e_{k}^{*} - \frac{1}{2} C \otimes \omega e'_{j}.
    \end{equation*}
    From these formulas one can see that
    \begin{equation*}
        R_{e_{i}}: E \rightarrow \mathbb{R}D, \ \ E'\rightarrow E\wedge E', \ \ E \wedge E' \rightarrow E', \ \ W \rightarrow W,
    \end{equation*}
    \begin{equation*}
        R_{e'_{j}}: E' \rightarrow \mathbb{R}D, \ \ E \rightarrow E \wedge E', \ \ E \wedge E' \rightarrow E, \ \ W \rightarrow W.
    \end{equation*}
    It implies that operators $R_{e_{i}}^{2}$, $R_{e'_{j}}^{2}$ can be considered as endomorphisms of the subspace $W$. So we have that
    \begin{equation*}
        \text{tr } R_{e_{i}}^{2} = \text{tr } R_{e_{i}}^{2}|_{W} = \frac{1}{4}\text{tr} \Big( C \otimes \omega e_{i}  \Big)^{2}
            = - \frac{1}{4} \text{tr} \Big( C^{2} \otimes \text{Id}_{\Sigma} \Big) = - \frac{d}{4} \text{tr}(C^{2}),
    \end{equation*}
    and, analogously, $\text{tr } R_{e'_{j}}^{2} = - \frac{d}{4} \text{tr}(C^{2})$. Substituting these values to the expressions
    (\ref{Ric_E_E'}) we obtain
    \begin{equation*}
        \text{Ric}(e_{i}, e_{i}) = \text{Ric}(e'_{j}, e'_{j}) = - \bigg( n + \frac{d}{2} \bigg) + \frac{d}{4} \text{tr}(C^{2}).
    \end{equation*}
    
    Finally, to compute the summand $\text{tr } R_{s}^{2}$ we note that
    \begin{equation*}
        R_{s} = L_{s} \circ \pi + \pi \circ \text{ad}_{s} = -\frac{1}{2} D \otimes s^{*} 
                + \frac{1}{2} \sum\limits_{i=1}^{p}\sum\limits_{j=1}^{q} (e_{i} \wedge e'_{j}) \otimes \big( (e_{i} e'_{j}) s \big)^{*}
    \end{equation*}
    \begin{equation*}
                + \frac{1}{2}\text{ad}_{s}|_{W}
                + \frac{1}{2} \sum\limits_{i=1}^{p} \Bigg( \Big( C \otimes \omega e_{i} \Big)  s \Bigg) \otimes e_{i}^{*}
                - \frac{1}{2}\sum\limits_{j=1}^{q} \Bigg( \Big( C \otimes \omega e'_{j} \Big) s \Bigg) \otimes {e'_{j}}^{*}.
    \end{equation*}
    From this one can see that
    \begin{equation*}
        R_{s}: E \rightarrow W, \ \ E' \rightarrow W, \ \ W \rightarrow \mathbb{R}D+E+E'+E \wedge E', 
    \end{equation*}
    which implies that only subspaces which give non-zero contribution to the computation of  $\text{tr}R_{s}^{2}$ are the subspaces $E$, $E'$, $W$.
    More precisely, for any $s \in W$ we have
    \begin{equation*}
        \text{tr } R_{s}^{2}|_{E} = \sum\limits_{i=1}^{p} g_{\hat{A}}(R_{s}^{2} e_{i}, e_{i}) 
            = \frac{1}{4} \sum\limits_{i=1}^{p} \langle [s, \big( C \otimes \omega e_{i} \big) s], e_{i} \rangle
            = - \frac{1}{4} \sum\limits_{i=1}^{p} \langle [ \big( C \otimes \omega e_{i} \big) s, s], e_{i} \rangle
    \end{equation*}
    \begin{equation*}
            = - \frac{1}{4} \sum\limits_{i=1}^{p} b_{A} \Big( \omega e_{i} \big( C \otimes \omega e_{i} \big) s, s \Big)
            = \frac{1}{4} \sum\limits_{i=1}^{p} b_{A} \Big( \big( C \otimes \text{Id}_{\Sigma} \big) s, s \Big)
            = \frac{p}{4} g_{\hat{A}} \Big( \big( C^{2} \otimes \text{Id}_{\Sigma} \big) s, s \Big)
    \end{equation*}
    \begin{equation*}
        \text{tr } R_{s}^{2}|_{E'} = \sum\limits_{j=1}^{q} g_{\hat{A}}(R_{s}^{2} e'_{j}, e'_{j}) 
            = \frac{1}{4} \sum\limits_{j=1}^{q} \langle [s, \big( C \otimes \omega e'_{j} \big) s], e'_{j} \rangle
            = - \frac{1}{4} \sum\limits_{j=1}^{q} \langle [ \big( C \otimes \omega e'_{j} \big) s, s], e'_{j} \rangle
    \end{equation*}
    \begin{equation*}
            = - \frac{1}{4} \sum\limits_{j=1}^{q} b_{A} \Big( \omega e'_{j} \big( C \otimes \omega e'_{j} \big) s, s \Big)
            = \frac{1}{4} \sum\limits_{j=1}^{q} b_{A} \Big( \big( C \otimes \text{Id}_{\Sigma} \big) s, s \Big)
            = \frac{q}{4} g_{\hat{A}} \Big( \big( C^{2} \otimes \text{Id}_{\Sigma} \big) s, s \Big)
    \end{equation*}
    To compute $\text{tr } R_{s}^{2}|_{W}$ we consider arbitrary basis $(s_{k})$ of $\Sigma$ which is orthonormal with respect to
    the bilinear form  $b$. Then $(\xi_{r} \otimes s_{k})$ is a basis of $W$. Also we denote the coefficients of the inverse matrix
    $\hat{A}^{-1}$ by $\hat{A}^{rm}$. we have the following sequence of equalities
    \begin{equation}
        \text{tr } R_{s}^{2}|_{W}
            = \sum\limits_{k = 1}^{d} \sum\limits_{r=1}^{l} \sum\limits_{m=1}^{l}
                    \hat{A}^{rm} g_{\hat{A}}(R_{s}^{2} \xi_{r} \otimes s_{k}, \xi_{m} \otimes s_{k})
    \end{equation}
    \begin{equation*}
            = \sum\limits_{k = 1}^{d} \sum\limits_{r=1}^{l} \sum\limits_{m=1}^{l}
                \hat{A}^{rm} g_{\hat{A}}(R_{s} [s, \xi_{r} \otimes s_{k}], \xi_{m} \otimes s_{k})
    \end{equation*}
    \begin{equation*}
            = \frac{1}{4} \sum\limits_{k = 1}^{d} \sum\limits_{r=1}^{l} \sum\limits_{m=1}^{l} \sum\limits_{i=1}^{p}
                    \hat{A}^{rm} g_{\hat{A}}(\langle [s, \xi_{r} \otimes s_{k}], e_{i} \rangle \Big( C \otimes \omega e_{i} \Big) s, \xi_{m} \otimes s_{k})
     \end{equation*}
     \begin{equation*}
                + \frac{1}{4} \sum\limits_{k = 1}^{d} \sum\limits_{r=1}^{l} \sum\limits_{m=1}^{l} \sum\limits_{j=1}^{q}
                    \hat{A}^{rm} g_{\hat{A}}(\langle [s, \xi_{r} \otimes s_{k}] , e'_{j} \rangle \Big( C \otimes \omega e'_{j} \Big) s, \xi_{m} \otimes s_{k})
    \end{equation*}
    \begin{equation*}
            = \frac{1}{4} \sum\limits_{k = 1}^{d} \sum\limits_{r=1}^{l} \sum\limits_{m=1}^{l} \sum\limits_{i=1}^{p}
                    \hat{A}^{rm} g_{\hat{A}}( b_{A}( \omega e_{i} s, \xi_{r} \otimes s_{k}) \Big( C \otimes \omega e_{i} \Big) s, \xi_{m} \otimes s_{k})
    \end{equation*}
    \begin{equation*}
                +\frac{1}{4} \sum\limits_{k = 1}^{d} \sum\limits_{r=1}^{l} \sum\limits_{m=1}^{l} \sum\limits_{j=1}^{q}
                    \hat{A}^{rm} g_{\hat{A}}(b_{A}( \omega e'_{j} s, \xi_{r} \otimes s_{k} ) \Big( C \otimes \omega e'_{j} \Big) s, \xi_{m} \otimes s_{k})
    \end{equation*}
    \begin{equation*}
            = \frac{1}{4} \sum\limits_{i=1}^{p} g_{\hat{A}}\Big(  \big( C \otimes \omega e_{i} \big) s,
                    \sum\limits_{k = 1}^{d} \sum\limits_{r=1}^{l} \sum\limits_{m=1}^{l} \hat{A}^{rm} b_{A}( \omega e_{i} s, \xi_{r} \otimes s_{k}) \xi_{m} \otimes s_{k} \Big)
    \end{equation*}
    \begin{equation*}
                +\frac{1}{4} \sum\limits_{j=1}^{q} g_{\hat{A}}\Big( \big( C \otimes \omega e'_{j} \big) s,
                    \sum\limits_{k = 1}^{d} \sum\limits_{r=1}^{l} \sum\limits_{m=1}^{l} \hat{A}^{rm} b_{A}( \omega e'_{j} s, \xi_{r} \otimes s_{k} ) \xi_{m} \otimes s_{k} \Big)
    \end{equation*}
    \begin{equation*}
            = \frac{1}{4} \sum\limits_{i=1}^{p} g_{\hat{A}}\Big(  \big( C \otimes \omega e_{i} \big) s,
                    \sum\limits_{k = 1}^{d} \sum\limits_{r=1}^{l} \sum\limits_{m=1}^{l} \hat{A}^{rm} g_{\hat{A}}\big( ( C \otimes \omega e_{i}) s, \xi_{r} \otimes s_{k} \big) \xi_{m} \otimes s_{k} \Big)
    \end{equation*}
    \begin{equation*}
                +\frac{1}{4} \sum\limits_{j=1}^{q} g_{\hat{A}}\Big( \big( C \otimes \omega e'_{j} \big) s,
                    \sum\limits_{k = 1}^{d} \sum\limits_{r=1}^{l} \sum\limits_{m=1}^{l} \hat{A}^{rm} g_{\hat{A}}\big( ( C \otimes \omega e'_{j}) s, \xi_{r} \otimes s_{k} \big) \xi_{m} \otimes s_{k} \Big)
    \end{equation*}
    \begin{equation*}
            = \frac{1}{4} \sum\limits_{i=1}^{p}
                    g_{\hat{A}}\big( (C \otimes \omega e_{i}) s, (C \otimes \omega e_{i}) s \big)
                + \frac{1}{4} \sum\limits_{j=1}^{q}
                    g_{\hat{A}}\big( (C \otimes \omega e'_{j}) s, (C \otimes \omega e'_{j}) s \big)
    \end{equation*}
    \begin{equation*}
            = \frac{1}{4} \sum\limits_{i=1}^{p}
                    g_{\hat{A}}\big( (C \otimes \text{Id}_{\Sigma}) s, (C \otimes \text{Id}_{\Sigma}) s \big)
                + \frac{1}{4} \sum\limits_{j=1}^{q}
                    g_{\hat{A}}\big( (C \otimes \text{Id}_{\Sigma}) s, (C \otimes \text{Id}_{\Sigma}) s \big)
    \end{equation*}
    \begin{equation*}
        = \frac{n}{4} g_{\hat{A}}\big( (C \otimes \text{Id}_{\Sigma}) s, (C \otimes \text{Id}_{\Sigma}) s \big)
        = \frac{n}{4} g_{\hat{A}}\big( (C^{2} \otimes \text{Id}_{\Sigma}) s, s \big).
    \end{equation*}
    Therefore, we have that    
    \begin{equation}
        \text{tr }R_{s}^{2} = \text{tr } R_{s}^{2}|_{E} + \text{tr } R_{s}^{2}|_{E'} + \text{tr } R_{s}^{2}|_{W}
            = \frac{n}{2} g_{\hat{A}} \Bigg( \Big( C^{2} \otimes \text{Id}_{\Sigma} \Big) s, s \Bigg).
    \end{equation}    
    Substituting this value for the trace to the expression (\ref{Ric_S}) we obtain that
    \begin{equation}
        \text{Ric}(s, s) = - \bigg( \frac{n}{2} + \frac{ld}{4} \bigg) g_{\hat{A}}(s, s)
                - \frac{n}{2} g_{\hat{A}} \Bigg( \Big( C^{2} \otimes \text{Id}_{\Sigma} \Big) s, s \Bigg).
    \end{equation}

    One of the Einstein equations, $\frac{\text{Ric}(s, s)}{g_{\hat{A}}(s, s)} = \frac{\text{Ric}(D, D)}{g_{\hat{A}}(D, D)} = - \big( n + \frac{ld}{4} \big)$,
    leads to the following condition on the endomorphism $C$
    \begin{equation*}
        g_{\hat{A}} \big( (C^{2} \otimes \text{Id}_{\Sigma}) s, s \big) = g_{\hat{A}}(s, s) \text{ for any } s \in W.
    \end{equation*}
    Since the endomorphism $C^{2} \otimes \text{Id}_{\Sigma}$ is symmetric with respect to $g_{\hat{A}}$ the condition above
    is equivalent to the equation $C^{2} = \text{Id}$. Also note that if this equation holds then $\text{tr }C^{2} = l$.
    Therefore, for this case we have
    \begin{equation*}\label{Ric_equality}
      \text{Ric}(D, D) = \text{Ric}(e_{i} \wedge e'_ {j}, e_{i} \wedge e'_{j}) = - \bigg( n + \frac{ld}{4} \bigg),
    \end{equation*}
    \begin{equation*}
        \text{Ric}(e_{i}, e_{i}) = \text{Ric}(e'_{j}, e'_{j}) = - \bigg( n + \frac{ld}{4} \bigg),
    \end{equation*}
    \begin{equation*}
        \text{Ric}(s, s) = - \bigg( n + \frac{ld}{4} \bigg) g_{\hat{A}}(s, s).
    \end{equation*}

    Since $\text{Ad}_{K}$-submodules $E \wedge E'$, $E$, $E'$ are irreducible the $\text{Ad}_{K}$-invariant bilinear forms
    $\lambda \cdot g_{\hat{A}}$ and $\text{Ric}$ coincide on each of them, but not only on the basis elements. Moreover,
    for $q > 1$ the submodules $\mathbb{R}D$, $E \wedge E'$, $E$, $E'$, $W$ are pairwise non-equivalent so they are
    orthogonal with respect to both $\lambda \cdot g_{\hat{A}}$ and $\text{Ric}$. For $q = 1$ we only need to check that
    $\mathbb{R} D \perp E'$ and $E \perp E \wedge E'$, but it can be done in the same way as in \cite[Th. 5]{AL}.
    \hfill$\Box$

\section{Isometries}\label{IsometriesSection}

    \begin{Theorem}\label{CongruenceIsometry}
        For any invertible matrix $P \in \text{\emph{GL}}(l)$ there exists an isometry
        $\Big( M_{p, q}(A), g_{M, \hat{A}} \Big) \simeq \Big( M_{p, q}(P^{T}AP), g_{M, P^{T}\hat{A}P}  \Big)$.
    \end{Theorem}
    \noindent
    \textit{Proof.} 
        Considered as a linear map $P: \mathbb{R}^{l} \longrightarrow \mathbb{R}^{l}$ transforms the bilinear form $\epsilon_{A}$
        to the bilinear form $\epsilon_{P^{T} A P}$. The linear map $P$ can be extended to the automorphism $P \otimes \text{Id}$
        of the module $W = \mathbb{R}^{l} \otimes \Sigma$. According to the construction, the automorphism $P \otimes \text{Id}$
        transforms the $V$-valued bilinear form $\Pi_{A}: \wedge^{2} W \longrightarrow V$ to $\Pi_{P^{T} A P}$. Therefore,
        we obtain the isomorphism of the extended Poincar\'{e} algebras $\mathfrak{p}_{p, q}(A) \simeq \mathfrak{p}_{p, q}(P^{T} A P)$
        which can be further extended to diffeomorphism of the homogeneous spaces $M_{p, q}(A) \simeq M_{p, q}(P^{T}AP)$.
        By the same reasons, the considered automorphism $P \otimes \text{Id}$ transforms the scalar product $g_{\hat{A}}$ to
        $g_{P^{T} \hat{A} P}$. Consequently, the constructed diffeomorphism is an isometry.
    \hfill$\Box$

    \begin{Theorem}\label{DiagonalizationTheorem}
        Einstein pseudo-Riemannian homogeneous space $\Big( M_{p, q}(A), g_{M, \hat{A}} \Big)$ is isometric to
        $\Big( M_{p, q}(A'), g_{M, \hat{A}'} \Big)$ where  $A', \hat{A}'$ are diagonal matrices with $\pm 1$ on the diagonal.
    \end{Theorem}
    \textit{Proof.}
        According to Theorem \ref{EinsteinCriterion}, the endomorphism $C$ is an involution so it can be diagonalized.
        From this we obtain that the bilinear forms $\epsilon_{A}$ and $\epsilon_{\hat{A}}$ on $\mathbb{R}^{l}$ can be
        simultaneously diagonalized by appropriate choice of basis (see \cite{Gre}). In other words, there exists an invertible
        matrix $P \in \text{GL}(l)$ such that $A' = P^{T}AP$, $\hat{A}' = P^{T}\hat{A}P$ are diagonal matrices. One can
        require that the diagonal of $A'$ consists of the elements $\pm 1$. Also note that we have the equality
        \begin{equation*}
            (\hat{A}'^{-1} A')^{2} = (P^{-1} \hat{A}^{-1} A P)^{2} = P^{-1} (\hat{A}^{-1} A)^{2} P = P^{-1} P = \text{Id}.
        \end{equation*}
        From this we obtain that $\hat{A}'^{2} = A'^{2}$ which implies that the diagonal of  $\hat{A}'$ consists of $\pm 1$ as well.
        Theorem \ref{CongruenceIsometry} gives the corresponding isometry.
    \hfill$\Box$
    
    From this theorem we obtain that one can associate with each Einstein homogeneous space $\Big( M_{p, q}(A), g_{M, \hat{A}} \Big)$
    the tuple of non-negative integers in the following way. Consider the basis of $\mathbb{R}^{l}$ in which the bilinear forms
    $\epsilon_{A}$, $\epsilon_{\hat{A}}$ are given by the diagonal matrices $A'$, $\hat{A}'$, respectively. Let $m_{\varepsilon, \eta} \in \mathbb{Z}_{\geq 0}$,
    where $\epsilon, \eta \in \{ +, - \}$, be the number of indicies $i \in \{1, ..., l \}$ such that $A'_{ii}$ has the sign $\varepsilon$,
    while $\hat{A}'_{ii}$ has the sign $\eta$. Note that any simultaneous congruent transformation of the pair of matrices
    $(A', \hat{A}')$ which respects the diagonal form does not change the numbers $m_{\varepsilon, \eta}$. Therefore,
    the 4-tuple $m := (m_{\varepsilon, \eta})$ is uniquely determined by the pair of matrices $(A, \hat{A})$.
    The immediate consequence of \cite[Th. 6]{C} is the following statement.
    \begin{Theorem}\label{TupleParametrization}
        For fixed $p \equiv 3 \pmod{4}$, $q \not\equiv 3 \pmod{4}$ the Einstein pseudo-Riemannian homogeneous space 
        $\Big( M_{p, q}(A), g_{M, \hat{A}} \Big)$ is completely determined by the tuple $m = (m_{\varepsilon, \eta})$ of four
        non-negative integers  up to isometry so we can denote $M(p, q, m) := \Big[ \Big( M_{p, q}(A), g_{M, \hat{A}} \Big) \Big]$.
        Moreover, if $m'_{\varepsilon, \eta} = m_{-\varepsilon, -\eta}$ then $M(p, q, m) = M(p, q, m')$. Also $M(p, q, m)$ is Riemannian
        if and only if $m_{\scriptscriptstyle{+, -}} = m_{\scriptscriptstyle{-, +}} = 0$.
    \end{Theorem}
    
\section{The case $p, q \equiv 3 \pmod{4}$}

    All methods and techniques considered above can be successfully applied to the case $p, q \equiv 3 \pmod{4}$ as well.
    The main difference is that for $p, q \equiv 3 \pmod{4}$ there exists exactly two irreducible non-equivalent
    $\text{Cl}^{0}(V)$-modules $\Sigma^{+}$ and $\Sigma^{-}$. These modules satisfy the relation $(\Sigma^{+})^{*} \simeq \Sigma^{-}$
    and have the same dimension $d = 4 \cdot N(p+q-3)$. Any $\text{Cl}^{0}(V)$-module $W$ can be represented as
    $W = W^{+} \oplus W^{-}$ where $W^{+} \simeq \mathbb{R}^{l^{+}} \otimes \Sigma^{+} $ and $W^{-} \simeq \mathbb{R}^{l^{-}} \otimes \Sigma^{-}$.
    By the same reasons as for $q \not\equiv 3 \pmod{4}$ there exists an unique (up to multiplication by positive scalar)
    $\hat{K}$-invariant symmetric positive-definite $\mathbb{R}$-valued bilinear form on each of the modules  $\Sigma^{+}$
    and $\Sigma^{-}$ which we denote by $b^{+}$ and $b^{-}$, respectively. Then any $\hat{K}$-invariant symmetric bilinear
    form on $W$ has the form $b_{A_{+}, A_{-}} = ( \epsilon_{A_{+}} \otimes b^{+} ) \oplus (\epsilon_{A_{-}} \otimes b^{-})$ where
    $\epsilon_{A_{+}}$, $\epsilon_{A_{-}}$ are some symmetric bilinear forms on $\mathbb{R}^{l^{+}}$ and $\mathbb{R}^{l^{-}}$,
    respectively, given by the matrices $A_{+}$, $A_{-}$ in the standard bases. Next, the natural modifications of the formulas
    (\ref{V-valuedForm}), (\ref{ExtendedPoincareAlgebra}), (\ref{Derivation}), (\ref{ExtensionOfEPA}) lead to
    the extended Pincar\'{e} algebra $\mathfrak{p}_{p, q}(A_{+}, A_{-})$, its $\mathbb{Z}_{2}$-graded extension $\mathfrak{g}_{p, q}(A_{+}, A_{-})$,
    and the corresponding Lie group $G_{p, q}(A_{+}, A_{-})$. In the similar way we can construct the homogeneous space
    $M_{p, q}(A_{+}, A_{-}) := G_{p, q}(A_{+}, A_{-})/K$ with the reductive decomposition of the form (\ref{ReductiveDecomposition}).
    We define ${Ad}_{K}$-invariant scalar product on the tangent space $\mathfrak{m}$ by analogy with (\ref{ScalarProduct})
    \begin{equation}
        g_{\hat{A}_{+}, \hat{A}_{-}} = g_{\mathbb{R}D} \oplus g_{E \wedge E'} \oplus g_{V} \oplus b_{\hat{A}_{+}, \hat{A}_{-}}.
    \end{equation}
    Finally, we obtain the pseudo-Riemannian homogeneous space \\ $\Big( M_{p,q}(A_{+}, A_{-}), g_{M, \hat{A}_{+}, \hat{A}_{-}} \Big)$.
    
    Denote the endomorphisms which transform the bilinear forms $\epsilon_{\hat{A}_{+}}$, $\epsilon_{\hat{A}_{-}}$ to
    the forms $\epsilon_{A_{+}}$, $\epsilon_{A_{-}}$ by $C_{+}$, $C_{-}$, respectively, i. e. 
    \begin{equation*}
        \epsilon_{\hat{A}_{+}}(C_{+} \xi, \xi' ) = \epsilon_{A_{+}}(\xi, \xi'), \ \ \
            \epsilon_{\hat{A}_{-}}(C_{-} \zeta, \zeta') = \epsilon_{A_{-}}(\zeta, \zeta').
    \end{equation*}
    Due to that $g_{\hat{A}_{+},\hat{A}_{-}}(W^{+}, W^{-}) = 0$ and $[W^{+}, W^{-}] = 0$ the analog of Theorem \ref{NomizuOperators} holds
    in which the summands
    \begin{equation*}
        \frac{1}{2} C \otimes \omega e_{i}, \ \ \frac{1}{2} C \otimes \omega e'_{j}, \ \
            \frac{1}{2}\sum\limits_{i=1}^{p} \Bigg( \Big( C \otimes \omega e_{i} \Big)  s \Bigg)  \otimes e_{i}^{*}, \ \
            \frac{1}{2}\sum\limits_{j=1}^{q} \Bigg( \Big( C \otimes \omega e'_{j} \Big) s  \Bigg) \otimes {e'_{j}}^{*}     
    \end{equation*}
    are replaced by
    \begin{equation*}
        \frac{1}{2} C_{+} \otimes \omega e_{i} + \frac{1}{2} C_{-} \otimes \omega e_{i}, \ \ \
        \frac{1}{2} C_{+} \otimes \omega e'_{j} + \frac{1}{2} C_{-} \otimes \omega e'_{j},
    \end{equation*}
    \begin{equation*}
        \frac{1}{2}\sum\limits_{i=1}^{p} \Bigg( \Big( C_{+} \otimes \omega e_{i} + C_{-} \otimes \omega e_{i} \Big)  s \Bigg)  \otimes e_{i}^{*},
        \frac{1}{2}\sum\limits_{j=1}^{q} \Bigg( \Big( C_{+} \otimes \omega e'_{j} + C_{-} \otimes \omega e'_{j} \Big) s  \Bigg) \otimes {e'_{j}}^{*},
    \end{equation*}
    respectively. Next, by the same reasons, the following modification of Theorem \ref{EinsteinCriterion} holds.
    \begin{Theorem}\label{EinsteinCriterion_q3}
        For $p, q \equiv 3 \pmod{4}$ the pseudo-Riemannian homogeneous space $\Big( M_{p, q}(A_{+}, A_{-}), g_{M, \hat{A}_{+}, \hat{A}_{-}} \Big)$
        is Einstein if and only if $C_{+}^{2} = \text{\emph{Id}}$ and $C_{-}^{2} = \text{\emph{Id}}$.
        The Einstein constant is $\lambda = -n - \frac{1}{4}(l^{+} + l^{-}) d$.
    \end{Theorem}
    The results of Section \ref{IsometriesSection} also can be transferred to the considered case. 
    \begin{Theorem}\label{CongruenceIsometry_q3}
        For $p, q \equiv 3 \pmod{4}$ for any two invertible matrices $P_{+} \in \text{\emph{GL}}(l^{+})$,
        $P_{-} \in \text{\emph{GL}}(l^{-})$ we have the following isometry:
       \begin{equation*}
           \Big( M_{p, q}(A_{+}, A_{-}), g_{M, \hat{A}_{+}, \hat{A}_{-}} \Big) \simeq
               \Big( M_{p, q}(P_{+}^{T}A_{+} P_{+}, P_{-}^{T} A_{-} P_{-}), g_{M, P_{+}^{T}\hat{A}_{+}P_{+}, P_{-}^{T}\hat{A}_{-}P_{-}}  \Big).
       \end{equation*}
    \end{Theorem}
    \begin{Theorem}\label{DiagonalizationTheorem_q3}
        For $p, q \equiv 3 \pmod{4}$ any Einstein pseudo-Riemannian homogeneous space $\Big( M_{p, q}(A_{+}, A_{-}), g_{M, \hat{A}_{+}, \hat{A}_{-}} \Big)$
        is isometric to \\ $\Big( M_{p, q}(A'_{+}, A'_{-}), g_{M, \hat{A}'_{+}, \hat{A}'_{-}} \Big)$ such that $A'_{+}, \hat{A}'_{+}, A'_{-}, \hat{A}'_{-}$
        are diagonal matrices with the elements $\pm 1$ on the diagonal.
    \end{Theorem}
    For the two pairs of matrices $(A_{+}, \hat{A}_{+})$, $(A_{-}, \hat{A}_{-})$ which determine Einstein homogeneous space
    $\Big( M_{p, q}(A_{+}, A_{-}), g_{M, \hat{A}_{+}, \hat{A}_{-}} \Big)$ we define by $m_{\varepsilon, \eta}^{\mu}$, where 
    $\varepsilon, \eta, \nu \in \{+, -\}$, the number of indicies $i \in \{1, ..., l^{\nu}\}$ such that ${(A^{\mu})}'_{ii}$ has
    the sign $\varepsilon$, while ${(\hat{A}^{\mu})}'_{ii}$ has the sign $\eta$. So to each homogeneous space of the considered
    form we put into correspondence the tuple of eight non-negative integers $m := \big( m_{\varepsilon, \eta}^{\mu} \big)$.
    The immediate consequence of \cite[Th. 7]{C} is the following statement.
    \begin{Theorem}\label{TupleParametrization_q3}
        For fixed $p, q \equiv 3 \pmod{4}$ the Einstein pseudo-Riemannian homogeneous space $\Big( M_{p, q}(A_{+}, A_{-}), g_{M, \hat{A}_{+}, \hat{A}_{-}} \Big)$
        is completely determined up to isometry by the 8-tuple $m = \big( m_{\varepsilon, \eta}^{\eta} \big)$, so we can denote
        $M(p, q, m) := \Bigg[ \Big( M_{p, q}(A_{+}, A_{-}), g_{M, \hat{A}_{+}, \hat{A}_{-}} \Big) \Bigg]$. Moreover, if $m' \in \Gamma(m)$,
        where $\Gamma \simeq \mathbb{Z}_{2} \times \mathbb{Z}_{2}$ is the group generated by two involutions $\alpha$ and
        $\beta$ such that $\alpha(m)_{\varepsilon, \eta}^{\mu} = m_{-\varepsilon, -\eta}^{\mu}$,
        $\beta(m)_{\varepsilon, \eta}^{\mu} = m_{\varepsilon, \eta}^{-\mu}$, then $M(p, q, m) = M(p, q, m')$.
        Also $M(p, q, m)$ is Riemannian if and only if $m_{\scriptscriptstyle{+, -}}^{\mu} = m_{\scriptscriptstyle{-, +}}^{\mu} = 0$
        for any $\mu \in \{ +, - \}$.
    \end{Theorem}

\end{document}